\documentclass[11pt]{article}
\usepackage{latexsym}
\usepackage{graphicx}
\usepackage{framed}
 \usepackage{threeparttable}

\begin{document}

\newcommand{\noi}{\noindent}
\newcommand{\nn}{\nonumber}
\newcommand{\bd}{\begin{displaymath}}
\newcommand{\ed}{\end{displaymath}}
\newcommand{\dfrac}[2]{\displaystyle\frac{#1}{#2}}
\newcommand{\R}{\rm I\kern-.19emR}
\newcommand{\bp}{\underline{\bf Proof}:\ }
\newcommand{\ep}{{\hfill $\Box$}\\ }
\newtheorem{1}{LEMMA}
\newtheorem{2}{THEOREM}
\newtheorem{3}[1]{COROLLARY}
\newtheorem{4}[1]{PROPOSITION}
\newtheorem{5}{REMARK}
\newtheorem{20}[1]{OBSERVATION}
\newtheorem{10}{DEFINITION}
\newtheorem{30}{RESULTS}
\newtheorem{40}[1]{CLAIM}
\newtheorem{50}[1]{ASSUMPTION}
\newtheorem{60}{EXAMPLE}
\newtheorem{70}{ALGORITHM}
\newcommand{\be}{\begin{equation}}
\newcommand{\ee}{\end{equation}}
\newcommand{\ba}{\begin{array}}
\newcommand{\ea}{\end{array}}
\newcommand{\bea}{\begin{eqnarray}}
\newcommand{\eea}{\end{eqnarray}}
\newcommand{\bqn}{\begin{eqnarray*}}
\newcommand{\eqn}{\end{eqnarray*}}

\newcommand{\e} { \ = \ }
\newcommand{\leqs}{ \ \leq \ }
\newcommand{\geqs}{ \ \geq \ }
\def\NN{\cal N}
 \def\C{{\rm\kern.24em
\vrule width.02em height1.4ex depth-.05ex \kern-.26em C}}
\def\theequation{\thesection.\arabic{equation}}
\def\reff#1{{\rm
(\ref{#1})}}
\def\E{{\mathbb{E}}\,}
\def\N{\mathbb{N}}
\def\P{{\mathbb{P}}}
\def\al{\alpha}
\newcommand{\eps}{\varepsilon}
\newcommand{\sgn}{\operatorname{sgn}}
\newcommand{\sign}{\operatorname{sign}}
\newcommand{\Vol}{\operatorname{Vol}}
\newcommand{\Var}{\operatorname{Var}}
\newcommand{\Cov}{\operatorname{Cov}}
\newcommand{\vol}{\operatorname{vol}}
\newcommand{\var}{\operatorname{var}}
\newcommand{\cov}{\operatorname{cov}}
\renewcommand{\Re}{\operatorname{Re}}
\renewcommand{\Im}{\operatorname{Im}}
\newcommand{\EE}{\mathbb E}
\newcommand{\CC}{\mathbb C}
\newcommand{\PP}{\mathbb P}
\newcommand{\QQ}{\mathbb Q}
\newcommand{\ZZ}{\mathbb Z}
\newcommand{\cA}{{\mathcal A}}
\newcommand{\cB}{{\mathcal B}}
\newcommand{\cC}{{\mathcal C}}
\newcommand{\cD}{{\mathcal D}}
\newcommand{\cE}{{\mathcal E}}
\newcommand{\cF}{{\mathcal F}}
\newcommand{\cG}{{\mathcal G}}
\newcommand{\cH}{{\mathcal H}}
\newcommand{\cI}{{\mathcal I}}
\newcommand{\cJ}{{\mathcal J}}
\newcommand{\cK}{{\mathcal K}}
\newcommand{\cL}{{\mathcal L}}
\newcommand{\cM}{{\mathcal M}}
\newcommand{\cN}{{\mathcal N}}
\newcommand{\cO}{{\mathcal O}}
\newcommand{\cP}{{\mathcal P}}
\newcommand{\cQ}{{\mathcal Q}}
\newcommand{\cR}{{\mathcal R}}
\newcommand{\cS}{{\mathcal S}}
\newcommand{\cT}{{\mathcal T}}
\newcommand{\cU}{{\mathcal U}}
\newcommand{\cV}{{\mathcal V}}
\newcommand{\cW}{{\mathcal W}}
\newcommand{\cX}{{\mathcal X}}
\newcommand{\cY}{{\mathcal Y}}
\newcommand{\cZ}{{\mathcal Z}}
\newcommand{\bx}{{\mathbf x}}
\newcommand{\by}{{\mathbf y}}
\newcommand{\bz}{{\mathbf z}}
\newcommand{\bba}{{\mathbf a}}
\newcommand{\bbb}{{\mathbf b}}
\newcommand{\bbc}{{\mathbf c}}


\title{An unconstrained optimization approach for finding real eigenvalues of even order symmetric tensors}
\author{
Lixing Han\thanks{
 Department of Mathematics, University of Michigan-Flint 
Flint, MI 48502, USA. {\it Email address}: {\tt lxhan@umflint.edu}
 }
}
\date{}
\maketitle

\begin{abstract} Let $n$ be a positive integer and $m$ be a  positive even integer. Let ${\mathcal A}$  be an $m^{th}$ order $n$-dimensional real weakly symmetric tensor and ${\mathcal B}$ be a real weakly symmetric positive definite tensor of the same size. $\lambda \in \R$ is called a ${\mathcal B}_r$-eigenvalue of ${\mathcal A}$ if ${\mathcal A} x^{m-1} = \lambda {\mathcal B} x^{m-1}$ for some 
$x \in \R^n \backslash \{0\}$. In this paper, we introduce two unconstrained optimization problems and obtain some variational characterizations for the minimum and maximum ${\mathcal B}_r$--eigenvalues of ${\mathcal A}$. Our results extend Auchmuty's unconstrained variational principles for eigenvalues of real symmetric matrices.  This unconstrained optimization approach can be used to find a Z-, H-, or D-eigenvalue of an even order weakly symmetric  tensor. We provide some numerical results to illustrate the effectiveness of this approach for finding a Z-eigenvalue and for determining the positive semidefiniteness of an even order symmetric tensor. 
\end{abstract}

\vskip 1cm
\ \\
{\bf Key words.} Weakly symmetric tensors,  tensor eigenvalues, positive semi-definiteness, unconstrained optimization.

\ \\
{\bf AMS subject classification (2010).} 65F15, 65K05, 15A69.

\section{Introduction}
\label{Intro}
\setcounter{equation}{0}

Since the pioneering works of Qi \cite{Qi} and Lim \cite{Lim}, the tensor eigenproblem has become an important part of numerical multilinear algebra.  In this paper, we consider the real eigenvalue problems for even order real symmetric tensors. Eigenvalues of symmetric tensors have found applications in several areas, including automatic control, statistical data analysis, higher order diffusion tensor imaging, and image authenticity verification, etc., see for example, \cite{Qi, QSW, QWWu, QYW, QYX}.  

Throughout this paper, we assume that  $\R$ is the real field, $n$ is a positive integer, and $m$ is a positive even integer. An  $m^{th}$-order $n$-dimensional real tensor   
$$ 
\cA = (A_{i_1i_2 \cdots i_m}) \in \R^{n \times n \times \cdots \times n}
$$  
 is called symmetric if its entries are invariant under any permutations of their indices \cite{KR,Qi}. A  tensor $\cA$ is called positive definite (positive semidefinite) if the multilinear form 
$$
 \cA x^m = \sum_{i_1, \cdots, i_m =1}^{n} A_{i_1i_2 \cdots i_m}x_{i_1} x_{i_2} \cdots x_{i_m} 
$$
is positive (nonnegative) for all $x \in \R^n \backslash \{ 0\}$.   The notation $\cA x^{m-1} $  denotes the vector in $\R^n$ whose  $i^{th}$ entry is
$$
(\cA x^{m-1})_i = \sum_{i_2, \cdots, i_m=1}^{n} A_{i i_2 \cdots i_m} x_{i_2}  \cdots x_{i_m}.
$$
 Following \cite{CPZ09}, $\cA$ is called weakly symmetric if  the gradient of $\cA x^m$
$$
 \nabla \left ( \cA x^m \right ) = m \cA x^{m-1}
$$
for all $x \in \R^n$. If $\cA$ is symmetric, then it is weakly symmetric \cite{CPZ09}.

Various definitions of real eigenpairs for tensors have been introduced in the literature, including H-eigenvalues \cite{Qi}, Z-eigenvalues \cite{Qi}, and D-eigenvalues  
\cite{QWWu}.   We use the following generalized eigenvalue definition, which includes the  H-, Z-, and D-eigenvalues as special cases.    
\begin{10}  
  Let $\cA$ and $\cB$  be $m^{th}$-order $n$-dimensional real weakly symmetric  tensors. Assume further that $\cB$ is positive definite.   If there exist a scalar $\lambda \in \R$ and a nonzero vector $ x  \in \R^n$ such that
\be
\label{evalue}
\cA x^{m-1} = \lambda \cB x^{m-1},
\ee
 then $\lambda$ is called a $\cB_{r}$--eigenvalue of $\cA$ and $x$ a $\cB_{r}$--eigenvector with respect to $\lambda$.  We denote the $\cB_{r}$--spectrum of $\cA$ by $$
\sigma_{\cB_{r}}({\cA}) = \{ \lambda : \lambda \ {\rm  is \ a \ \cB_{r}{\rm -}eigenvalue \ of }\ \cA \}.
$$ 
\end{10}

\begin{5} 
\label{remark1}
{\rm
This definition was first introduced by Chang, Pearson, and Zhang \cite{CPZ09} in a somewhat more general setting. 
\begin{itemize}
\item If $\cB =\cI = (\delta_{i_1 i_2 \cdots i_m})$, the unit tensor, then $\cB$ is weakly symmetric  positive definite. Moreover,  $\cB x^m = \|x\|_{m}^m = x_{1}^m + x_{2}^m + \cdots + x_{n}^m$,
$\cB x^{m-1}= m[x_{1}^{m-1}, x_{2}^{m-1}, \cdots, x_{n}^{m-1}]^T $, and the $\cB_r$--eigenvalues are H-eigenvalues. 
\item If $\cB = I_{n}^{m/2}$, the tensor product of $m/2$ copies of the unit matrix $I_n \in \R^{n\times n} $, then $\cB$ is weakly symmetric  positive definite. Moreover, 
$\cB x^m = (x^Tx)^{m/2}=(x_{1}^2 + x_{2}^2 + \cdots + x_{n}^2)^{m/2}   $, $\cB x^{m-1} = m (x^Tx)^{\frac{m-2}{2}} x $, and  the $\cB_r$--eigenvalues are Z-eigenvalues.  
\item If $\cB = D^{m/2}$, the tensor product of $m/2$ copies of the symmetric positive definite  matrix $D \in \R^{n\times n} $, 
then  $\cB$ is weakly symmetric  positive definite. Moreover, 
$\cB x^m =  (x^TDx)^{m/2}$,  $\cB x^{m-1} = m (x^TDx)^{\frac{m-2}{2}} Dx $, and
the $\cB_r$--eigenvalues are D-eigenvalues. 
\end{itemize}
}
\end{5}

Calculation of all eigenvalues of a high order ($m>2$) tensor is difficult, unless $m$ and $n$ are small \cite{Qi}. In certain circumstances, however, one only needs to compute the largest or smallest eigenvalue of a tensor. For instance, the smallest H-eigenvalue or Z-eigenvalue of an even order symmetric tensor $\cA$ can be used to determine the positive definiteness/semidefiniteness of $\cA$ \cite{Qi}. For a nonnegative tensor, the Perron-Frobenius theory asserts that its largest H-eigenvalue is its spectral radius \cite{CPZ08, FGH}.

Recently, Kolda and Mayo \cite{KM} have extended the high order power method for symmetric tensor eigenproblems of Kofidis and Regalia \cite{KR} by introducing a shift parameter $\alpha$ to compute Z-eigenvalues of symmetric tensors.  With a suitable choice of $\alpha$, the resulting method,  SSHOPM, converges to a Z-eigenvalue of the tensor when applied to a symmetric tensor.  The found Z-eigenvalue is not necessarily the largest or smallest Z-eigenvalue. The rate of convergence of the SSHOPM method is  linear \cite{KM}.

An alternative approach for computing the  eigenvalues of a symmetric tensor is to solve the constrained optimization problem \cite{Lim}
\be
\label{minconopt}
\min \cA x^m \ \ {\rm s.t.} \ \ \cB x^m =1,   
\ee  
or
\be
\label{maxconopt}
\max \cA x^m \ \ {\rm s.t.} \ \ \cB x^m =1.
\ee
The Karush-Kuhn-Tucker points of Problem \reff{minconopt} or \reff{maxconopt} give $\cB_r$--eigenvalues and $\cB_r$--eigenvectors of $\cA$. If we are interested in obtaining one eigenvalue, then   these problems can be solved using a local constrained optimization solver \cite{NW}. Note that in each problem, the objective function and the constraint function are both polynomials. Therefore, a global polynomial optimization method  can be used,  if we are interested in finding the largest or smallest $\cB_r$--eigenvalue.

A more attractive approach for computing  eigenvalues of even order symmetric tensors,  however, is to use unconstrained optimization. This is motivated by the works of Auchmuty \cite{Au89, Au91}, in which he proposed some unconstrained variational principles for generalized symmetric matrix eigenvalue problems. In particular, he \cite{Au91} considered the unconstrained optimization problems
\be
\label{Auobj1}
\min_{x \in \R^n} g_1(x)= \frac{1}{4} (x^T B x)^2 + \frac{1}{2} x^T A x,  
\ee
and
\be
\label{Auobj2}
\min_{x \in \R^n} g_2(x)= \frac{1}{4} (x^T B x)^2 - \frac{1}{2}  x^TAx,  
\ee
where $A \in \R^{n \times n}$ is a symmetric matrix and $B \in \R^{n \times n}$ is a symmetric positive definite matrix. He proved that Problem \reff{Auobj1} can be used to find the smallest generalized $B$-eigenvalue of $A$ and Problem \reff{Auobj2} can be used to find the largest generalized $B$--eigenvalue of $A$.
In this paper, we will extend Auchmuty's unconstrained variational principles for symmetric matrix eigenproblems \cite{Au91}  to even order weakly symmetric tensors. 

The rest of this paper is organized as follows.  In Section 2, we introduce some preliminary results that will be used to establish the main results in Section 3.  In Section 3, we introduce two unconstrained optimization problems and obtain some variational characterizations for the
minimum and maximum $\cB_r$--eigenvalues of $\cA$. In Section 4, we give some numerical results. Some final remarks are given in Section 5.

\section{Preliminaries} 
\label{prelim}
\setcounter{equation}{0}

We start with the existence of $\cB_r$--eigenvalues of $\cA$.  The existence of H-eigenvalues and Z-eigenvalues of an even order symmetric tensor was first studied by Qi \cite{Qi}.  In \cite{CPZ09},  Chang, Pearson, and Zhang proved the existence of at least $n$ $\cB_r$--eigenvalues when $\cA$ is weakly symmetric and $\cB$ is weakly symmetric positive definite, which is summarized in the following   
\begin{2} (\cite{CPZ09})
\label{eignumber}
Assume that $\cA$ and $\cB$ are $m^{th}$-order n-dimensional real weakly symmetric tensors and $\cB$ is positive definite. Then $\cA$ has at least $n$  $\cB_r$--eigenvalues, with $n$ distinct pairs of $\cB_r$--eigenvectors. 
\end{2}

In \cite{Qi}, Qi  proved the existence of the maximum and minimum H-eigenvalues and Z-eigenvalues.  Using a similar argument, we can prove   
\begin{2}
\label{eigminmax}
Assume that $\cA$ and $\cB$ are $m^{th}$-order n-dimensional real weakly symmetric tensors and $\cB$ is positive definite. Then $\sigma_{\cB_r}(\cA)$ is not empty. Furthermore, 
there exist $ \lambda_{\min} \in \sigma_{\cB_r}(\cA)$ and $ \lambda_{\max} \in \sigma_{\cB_r}(\cA)$ such that
$$
  -\infty < \lambda_{\min} \leq \lambda \leq \lambda_{\max} < \infty, \ \ \forall \ \lambda \in \sigma_{\cB_r} (\cA). 
$$
\end{2}
\bp
Since $\cB$ is positive definite, the set $\{x \in \R^n: \cB x^m =1\}$  is compact \cite{CPZ09}. We also notice that function $\cA x^m$ is continuous. Thus, 
 the constrained optimization problem \reff{minconopt}
has a global minimizer $\underline{x}$ and  the constrained optimization problem \reff{maxconopt}
has a global maximizer $\bar{x}$. 

At the global minimizer $\underline{x}$ of problem \reff{minconopt}, there is a scalar $\underline{\lambda}$ such that the KKT conditions 
\be
\label{kkt1}
 m \cA \underline{x}^{m-1} = \underline{\lambda} m \cB \underline{x}^{m-1} 
\ee
hold. Clearly, $\underline{\lambda} \in \sigma_{\cB_r}(\cA)$. The inner product of \reff{kkt1} with $\underline{x}$ gives 
$$
\cA \underline{x}^m = \underline{\lambda}\cB \underline{x}^m =   \underline{\lambda}.
$$
Since $\underline{x}$ is a global minimizer of problem \reff{minconopt}, 
$$
\underline{\lambda} \leq \lambda, \ \ \ \forall \lambda \in \sigma_{\cB_r} (\cA). 
$$ 
Therefore, we can set $\lambda_{\min}=\underline{\lambda}$. Similarly, we can establish the existence of $\lambda_{\max}$ by using the global maximizer $\bar{x}$.  
\ep

We next consider a property of weakly symmetric positive definite tensors, which is similar to a  property for symmetric positive definite matrices. 
\begin{2}
\label{pdtensor}
Assume that $\cB$ is an $m^{th}$-order n-dimensional weakly symmetric positive definite tensor. Let $\mu>0$ be the smallest H-eigenvalue of $\cB$. Then 
\be
\label{ineq}
\cB x^m \geq \mu \|x\|_{m}^m, \ \ \forall x \in \R^n, 
\ee
where $\|x\|_{m}$ is the $m$-norm of $x$. 
\end{2}
\bp
When $x=0$, \reff{ineq} obviously holds. According to Theorem \ref{eigminmax},  $\mu$ is the global minimum value of
$$
  \min \cB x^m, \ \ \ {\rm s.t.} \ \ \|x\|_{m}^m=1. 
$$
For any $x \in \R^n \backslash \{0\}$, we have
$$
\cB \left ( \frac{x}{\|x\|_{m}} \right )^m \geq \mu. 
$$ 
This implies \reff{ineq}. 
\ep

Finally we recall that a continuous function $f: \R^n \to \R$ is {\it coercive} if 
$$
\lim_{\|x\| \to \infty} f(x) = + \infty.
$$
A nice feature of coercive functions is summarized in the following
\begin{2} (\cite{PSU})
\label{coerciveThm}
Let $f: \R^n \to \R$ be continuous. If $f$ is coercive, then $f$ has at least one global minimizer.   If, in addition, the first partial derivatives exist on $\R^n$, then $f$ attains its global minimizers at its critical points. 
\end{2}

\section{Unconstrained variational principles for  the minimal and maximal $\cB_r$ eigenvalues}
\label{Variational}
\setcounter{equation}{0}

We now generalize the unconstrained variational principles of Auchmuty \cite{Au91} to even order weakly symmetric tensors. We first consider the unconstrained optimization problem
\be
\label{obj1}
\min f_1(x)= \frac{1}{2m} (\cB x^m)^2 + \frac{1}{m} \cA x^m.  
\ee
When $\cA$ and $\cB$ are weakly symmetric, the gradient of the objective function $f_1$ is
\be
\label{gradient1}
\nabla f_1(x)=   (\cB x^m) \cB x^{m-1}     + \cA x^{m-1} .
\ee
The following theorem summarizes the properties of function $f_1$. 
\begin{2}
\label{mainthm1}
 Assume that $\cA$ and $\cB$ are $m^{th}$-order $n$-dimensional real weakly symmetric tensors and $\cB$ is positive definite. Let $\lambda_{\min}$ be the smallest $\cB_r$--eigenvalue of $\cA$. Then \\
(a) $f_1$ is coercive on $\R^n$. \\
(b) The  critical points of $f_1$ are  \\
\mbox{} \ \ \ (i) $x=0$; and  \\
\mbox{} \ \ \  (ii) any $\cB_r$--eigenvector $x$ of $\cA$ associated with a $\cB_r$--eigenvalue $\lambda < 0$ of $\cA$ satisfying $\cB x^m = - \lambda $.  \\
(c) If $\lambda_{\min} <0$, then $f_1$ attains its global minimal value  
$$
\min f_1(x) = -\frac{1}{2m} \lambda_{\min}^2
$$ 
at any $\cB_r$--eigenvector associated with the $\cB_r$--eigenvalue $\lambda_{\min}$ satisfying $\cB x^m =  - \lambda_{\min}$ . \\
(d) If $\lambda_{\min} \geq 0$, then $x=0$ is the unique critical point of $f_1$ and the unique global minimizer of $f_1$ on $\R^n$.  
\end{2}

\bp
(a) Since $\cB$ is weakly symmetric positive definite, Theorem \ref{pdtensor} asserts that
$$
 \cB x^m \geq \mu \|x\|_{m}^m, 
$$
where $\mu>0$ is the smallest H-eigenvalue of $\cB$.  This implies 
$$
f_1(x) \geq \frac{\mu^2}{2m} \|x\|_{m}^{2m} + \frac{1}{m} \cA x^m \to \infty  
$$
as $\|x\| \to \infty$. Thus, $f_1$ is coercive on $\R^n$. \\
(b) At a critical point of $f_1$, its gradient  $\nabla f_1(x)=0$, that is,  
\be
\label{eigeq1}
\cA x^{m-1} = -  (\cB x^m) \cB x^{m-1}. 
\ee
Clearly, $x=0$ is a critical point of $f_1$ as $\nabla f_1(0)=0$. Moreover, if $\lambda <0$ is a $\cB_r$--eigenvalue of $\cA$, then
$$
 \cA x^{m-1}= \lambda \cB x^{m-1}. 
$$
If $x \in \R^n \backslash \{ 0 \}$ is a  $\cB_r$--eigenvector associated with this $\lambda$ and satisfies $\cB x^m = - \lambda$, then it is a critical point of $f_1$.  \\
(c) From Theorem \ref{eigminmax}, $\lambda \geq \lambda_{\min}, \ \forall \lambda \in \sigma_{\cB_r} (\cA) $.  
At the critical point $x  \in \R^n\backslash \{0\}$ that is a  $\cB_r$--eigenvector associated with a $\cB_r$--eigenvalue $\lambda<0$ and satisfies $\cB x^m = - \lambda$, $ \cA x^m = - \lambda^2$. Moreover,
$$
f_1(x)= \frac{1}{2m} \lambda^2 - \frac{1}{m} \lambda^2 = - \frac{1}{2m} \lambda^2 \geq - \frac{1}{2m} \lambda_{\min}^2,  
$$  
since $0> \lambda \geq \lambda_{\min} $. According to Theorem \ref{coerciveThm} and part (b), $f_1$ attains the global minimum value $ \displaystyle{- \frac{1}{2m} \lambda_{\min}^2 }$ at  any $\cB_r$--eigenvector associated with the $\cB_r$--eigenvalue $\lambda_{\min}$ satisfying $\cB x^m =  - \lambda_{\min}$. \\
(d) $\lambda_{\min} \geq 0$ implies that $\lambda \geq 0$ for any $\lambda \in \sigma_{\cB_r}(\cA)$. Thus, $\cB x^m = - \lambda $ does not hold  for any $\cB_r$--eigenvector $x$ of $\cA$ associated with a $\cB_r$--eigenvalue $\lambda $ of $\cA$, as  $\cB x^m >0 $ for any $x \in \R^n \backslash \{0\}$ by the positive definiteness of $\cB$. Hence, $x=0$ is the unique critical point of $f_1$. It is also the unique 
global minimizer of $f_1$ according to Theorem \ref{coerciveThm}.     
\ep

We next consider the unconstrained optimization problem
\be
\label{obj2}
\min f_2(x)= \frac{1}{2m} (\cB x^m)^2 - \frac{1}{m} \cA x^m.  
\ee
When $\cA$ and $\cB$ are weakly symmetric, the gradient of the objective function $f_2$ is
\be
\label{gradient2}
\nabla f_2(x)=   (\cB x^m) \cB x^{m-1}     - \cA x^{m-1} .
\ee
Using a similar argument in the proof of the properties of $f_1$, we can prove the following properties about $f_2$.  
\begin{2}
\label{mainthm2}
 Assume that $\cA$ and $\cB$ are $m^{th}$-order n-dimensional real weakly symmetric tensors and $\cB$ is positive definite. Let $\lambda_{\max}$ be the largest $\cB_r$-eigenvalue of $\cA$. Then \\
(a) $f_2$ is coercive on $\R^n$. \\
(b) The  critical points of $f_2$ are at \\
\mbox{} \ \ \ \ (i) $x=0$; and  \\
\mbox{} \ \ \  (ii) any $\cB_r$--eigenvector $x$ of $\cA$ associated with a $\cB_r$--eigenvalue $\lambda > 0$ of $\cA$ satisfying $\cB x^m = \lambda$.  \\
(c) If $\lambda_{\max} >0$, then $f_2$ attains its global minimal value  
$$
\min f_2(x) = -\frac{1}{2m} \lambda_{\max}^2
$$ 
at any $\cB$-eigenvector associated with the $\cB_r$--eigenvalue $\lambda_{\max}$ satisfying $\cB x^m = \lambda_{\max}$. \\
(d) If $\lambda_{\max} \leq 0$, then $x=0$ is the unique critical point of $f_2$. Moreover, it is the unique global minimizer of $f_2$ on $\R^n$.  
\end{2}

Note that the functions $f_1$ and $f_2$ are polynomials of degree $2m$. A global polynomial optimization solver such as GloptiPoly3 \cite{HLL} can be used to find the smallest $\cB_r$--eigenvalue of $\cA$ and the largest $\cB_r$--eigenvalue of $\cA$ by solving Problem \reff{obj1} and Problem \reff{obj2} respectively, provided that $\lambda_{\min}<0$ and $\lambda_{\max} >0$. 

If $\lambda_{\min} \geq 0$, however, solving Problem \reff{obj1} does not result in the smallest $\lambda_{\min}$.  In this case, we can solve the shifted problem
\be
\label{shifted1}
\min_{x \in \R^n} s_1(x,t) = \frac{1}{2m} (\cB x^m)^2 + \frac{1}{m} (\cA + t \cB) x^m,  
\ee
using a  a suitable parameter $t<0$. Specifically, if  $t < - \lambda_{\min} $, then the global minimum value of Problem \reff{shifted1} is $\displaystyle{-\frac{1}{2m}(\lambda_{\min} +t)^2}$. Thus, $\lambda_{\min}$ can be obtained by finding the global minimum of Problem \reff{shifted1}. Summarizing the above discussions, we have the following theorem for Problem \reff{shifted1}. 
\begin{2}
\label{mainthm3}
Assume that $\cA$ and $\cB$ are $m^{th}$-order $n$-dimensional real weakly symmetric tensors and $\cB$ is positive definite. Let $\lambda_{\min}$ be the smallest $\cB_r$--eigenvalue of $\cA$. \\
(a) If $t+ \lambda_{\min} <0$, then the objective function $s_1$ defined in \reff{shifted1} attains its global minimal value  
$$
s_{1}^{*}=\min s_1(x) = -\frac{1}{2m} (t+\lambda_{\min})^2
$$ 
at any $\cB_r$--eigenvector of $\cA$ associated with the $\cB_r$--eigenvalue $\lambda_{\min}$ satisfying $\cB x^m =  - (t+\lambda_{\min})$. 
Moreover, we have 
$$
\lambda_{\min} = - \sqrt{-2ms_{1}^{*}}-t. 
$$ 
(b) If $t+ \lambda_{\min} \geq 0$, then $x=0$ is the unique critical point of $s_1$ and the unique global minimizer of $s_1$ on $\R^n$.  
\end{2} 

When $\lambda_{\max} \leq 0$, we can similarly solve the shifted problem
\be
\label{shifted2}
\min_{x \in \R^n} s_2(x,t) = \frac{1}{2m} (\cB x^m)^2 - \frac{1}{m} (\cA + t \cB) x^m,
\ee 
by using a suitable parameter $t >- \lambda_{\max} \geq 0$ to find $\lambda_{\max}$. We have the following theorem for Problem \reff{shifted2}.
\begin{2}
\label{mainthm4}
Assume that $\cA$ and $\cB$ are $m^{th}$-order $n$-dimensional real weakly symmetric tensors and $\cB$ is positive definite. Let $\lambda_{\max}$ be the largest $\cB_r$--eigenvalue of $\cA$. \\
(a) If $t+ \lambda_{\max} > 0$, then the objective function $s_2$ defined in \reff{shifted2} attains its global minimal value  
$$
s_{2}^{*} = \min s_2(x) = -\frac{1}{2m} (t+\lambda_{\max})^2
$$ 
at any $\cB_r$--eigenvector of $\cA$ associated with the $\cB_r$--eigenvalue $\lambda_{\max}$ satisfying $\cB x^m =  (t+\lambda_{\max})$. Moreover, we have
$$
\lambda_{\max} =  \sqrt{-2ms_{2}^{*}}-t. 
$$ 
(b) If $t+\lambda_{\max} \leq 0$, then $x=0$ is the unique critical point of $s_2$. Moreover, it is the unique global minimizer of $s_2$ on $\R^n$. 
\end{2} 

Problems \reff{obj1} and \reff{shifted1} can be used to determine whether an even order symmetric tensor $\cA$ is positive semidefinite or not.  Take $\cB=\cI \ {\rm or} \ \cB=I_{n}^{m/2}$. If the global minimum value of $f_1$ equals $0$, then $\cA$ is positive semidefinite (or definite); otherwise, it is not. Assume that we have been able to determine that $\cA$ is positive semidefinite.  To further determine whether $\cA$ is positive definite or semidefinite, we can solve \reff{shifted1} using $t=-1$. If the global minimum of $s_1$ is  $-\frac{1}{2m}$, then $\cA$ is only positive semidefinite;  otherwise $\cA$ is positive definite. 

Local unconstrained optimization methods can be used to solve Problems \reff{obj1} and \reff{obj2}. These methods do not guarantee finding a global minimum. However, they converge to a critical point (see for example, \cite{NW}). According to Theorems \ref{mainthm1} and \ref{mainthm2}, the found nonzero critical point corresponds to a $\cB_r$--eigenvalue of $\cA$. Therefore, local optimization solvers have the ability to find other eigenvalues besides the extreme ones. Moreover, if solving Problem \reff{obj1} with $\cB=I_{n}^{m/2}$ or $\cB=\cI$ results in a nonzero critical point, then it corresponds to a negative Z-eigenvalue or H-eigenvalue. This implies that local unconstrained optimization solvers can be used to solve Problem \reff{obj1} or Problem \reff{shifted1} to determine if $\cA$ is positive semidefinite. Finally, a local unconstrained optimization method such as the BFGS method has a fast rate of convergence - which is superlinear.

\section{Numerical results} 
\label{Num}
\setcounter{equation}{0}

In this section, we present some numerical results to illustrate the effectiveness of using the unconstrained variational principles for finding real eigenvalues of even order symmetric tensors.  The experiments were done on a laptop computer with an i3-2357M CPU @1.30GHz and a 4GB RAM running Windows 7, using MATLAB7.8.0 \cite{Mathworks}, the  MATLAB Optimization Toolbox \cite{Mathworks}, and the Tensor Toolbox \cite{BK}. We did two groups of experiments: First, comparing the new approach with the SSHOPM method and the constrained optimization approach. Second, testing the ability of the new approach to determine positive semidefiniteness of even order symmetric tensors.

\subsection{Effectiveness of finding a Z-eigenvalue}

In our first group of experiments, we tested the new approach on finding  Z-eigenvalues ($\cB=I_{n}^{m/2}$) in order to compare it with the SSHOPM method (\cite{KM}). We will focus on solving Problem \reff{obj1}. The numerical behavior of solving Problem \reff{obj2} is similar. 
When $\cB=I_{n}^{m/2}$, the unconstrained variational principle \reff{obj1} becomes  
\be
\label{Zfun1}
\min f_{1}^{Z}(x)= \frac{1}{2m} (x^Tx)^{m} +\frac{1}{m} \cA x^m.
\ee 
The gradient of the corresponding objective is 
\be 
\nabla f_{1}^Z(x) = (x^Tx)^{m-1}  x + \cA x^{m-1}.
\ee 
We tested the symmetric $4^{th}$ order tensors defined in the following examples: 

\begin{60} {\rm The $4^{th}$ order n-dimensional tensor $\cA$ is defined by
\label{example1}
\bea
 \cA(i,j,k,l) =  \left \{  \begin{array}{cl}
     - 0.9, & {\rm if} \ i=j=k=l; \\
\ \\
     0.1,  & {\rm otherwise}.
       \end{array}
\right . 1 \leq i,j,k,l \leq n
\eea
}
\end{60}

\begin{60}
\label{example2}
{\rm
The $4^{th}$ order n-dimensional symmetric tensor $\cA$ is generated as follows: First randomly generate tensor 
$ \cT ={\tt randn}(n,n,n,n)$, then  use the {\tt symmetrize} function in the Matlab Tensor Toolbox \cite{BK} to symmetrize $\cT$ and obtain   
$\cA= {\tt symmetrize} (\cT)$.
}
\end{60}

\begin{60}
\label{example3}
{\rm
The $4^{th}$ order n-dimensional symmetric tensor $\cA$ is generated as follows: First randomly generate tensor 
$ \cY ={\tt randn}(n,n,n,n)$; then create tensor $\cZ$ by setting $\cZ(i,j,k,l)= \frac{1}{\cY(i,j,k,l)} $; and finally  use the {\tt symmetrize} function in the Matlab Tensor Toolbox \cite{BK} to symmetrize $\cZ$ and obtain   
$\cA= {\tt symmetrize} (\cZ)$.
}
\end{60}

Since sometimes we are interested in finding the extreme eigenvalues of a tensor, we used the global polynomial optimization solver {\tt GloptiPoly3} of Henrion, Lasserre, and L\"{o}fberg \cite{HLL} to solve Problem \reff{Zfun1} for some $4^{th}$ order symmetric tensors in Examples 1--3.  We observed that {\tt GloptiPoly3} was able to solve \reff{Zfun1} when $n \leq 7$. When $n \geq 8$, it was unable to solve \reff{Zfun1} due to its memory requirement exceeding the capacity of the laptop computer we used.

From now on in this subsection we shall focus on solving \reff{Zfun1} using a local optimization method. Specifically,  we used the local optimization solver {\tt fminunc} (which uses a line search BFGS method) from the Matlab Optimization Toolbox \cite{Mathworks} to solve Problem \reff{Zfun1}, with its default settings except for the following: 
\be
\label{optset1}
{\tt GradObj:on, \ LargeScale:off, \ TolX = TolFun=10^{-12}, \ MaxIter=1000.}
\ee
We tested this approach and compared it with two other approaches on some tensors from Examples 1--3. 

\subsubsection{Comparison with the constrained variational principle \reff{minconopt}}

We  tested and compared  the unconstrained variational principle \reff{Zfun1} with the constrained variational approach \reff{minconopt} (using $\cB = I_{n}^{m/2}$) for finding Z-eigenvalues of some $4^{th}$ order n-dimensional tensors given in Examples 1 and 2. For the constrained variational principle approach, we used the {\tt fmincon} function from the Matlab Optimization 
Toolbox \cite{Mathworks} to solve Problem \reff{minconopt}.  We used  the default settings of {\tt fmincon} except for the following: 
\be
\label{optset2}
{\tt GradObj:on, \ GradConstr:on, \ TolX=TolFun=10^{-12}, \ MaxIter=1000.}
\ee

Our numerical experiments have shown that Problem \reff{minconopt} can be a surprisingly difficult problem for {\tt fmincon} when $m=4$. Take  the tensor $\cA$ in Example \ref{example1} with $n=4$ as an example. We  ran both {\tt fminunc} and {\tt fmincon} on this tensor, using  randomly generated initial vectors $x_0={\tt randn}(4,1)$ 100 times and using normalized randomly generated initial vectors $x_0 = y_0/\|y_0\|_2$ 100 times, where $y_0= {\tt randn}(4,1)$.  We observed that
\begin{itemize}
 \item Solving \reff{Zfun1}  via {\tt fminunc}:  In all of 200 runs, this method successfully found the Z-eigenvalue $\lambda=-0.9345$ of $\cA$.  
 \item Solving \reff{minconopt} with $\cB = I_{n}^{m/2}$  via {\tt fmincon}: (a) In the 100 runs using randomly generated initial vectors, this method successfully found the Z-eigenvalue  $\lambda=-0.9345$ of $\cA$ 11 times and it failed to find a Z-eigenvalue of $\cA$ in 89 runs. (b) In the 100 runs using normalized randomly generated initial vectors, it successfully found the Z-eigenvalue  $\lambda=-0.9345$ of $\cA$ 43 times and it failed to find a Z-eigenvalue of $\cA$ in 57 runs. The failures in both case (a) and case (b) were due to the divergence of {\tt fmincon}. {\tt fmincon} is based on a Sequential Quadratic Programming (SQP) method for constrained nonlinear optimization. It is generally robust. However,  \reff{minconopt} seems rather difficult for {\tt fmincon}. This indicates that new algorithms need to be developed to solve the constrained problem \reff{minconopt} directly. 
\end{itemize}

\subsubsection{Comparison with the SSHOPM method}

We now compare the unconstrained variational principle \reff{Zfun1} with the SSHOPM method of Kolda and Mayo \cite{KM} for finding Z-eigenvalues of some $4^{th}$ order n-dimensional tensors. The SSHOPM method is implemented as the {\tt sshopm} function in the Matlab Tensor Toolbox \cite{BK}. We  used the default settings of  {\tt sshopm} except for the following:
\be
\label{optset3}
{\tt  Tol=10^{-12}, \ MaxIts=5000.}
\ee

The iterates $x^{(k)}$ generated by {\tt sshopm} keeps the norm $\| x^{(k)} \|_2 =1$. When {\tt fminunc} converges to a nonzero critical point $\tilde{x}$ of Problem \reff{Zfun1} corresponding to a Z-eigenvalue $\lambda \leq 0$, $\|\tilde{x}\|_{2}^m = -\lambda$. To  have a fair comparison of the two methods, we first normalized the nonzero vector $\tilde{x}$ obtained by {\tt fminunc} at termination so that $\hat{x}= \tilde{x}/\|\tilde{x}\|_2$.  We define the error term by
\be
\label{success}
\hat{\epsilon} =  \| \cA \hat{x}^{3} - \lambda \hat{x} \|_2,
\ee
where $\hat{x}$ is either the vector obtained by {\tt sshopm} at termination or the normalized vector when solving \reff{Zfun1} via {\tt fminunc} at termination.

The SSHOPM method without shift (i.e., the original SHOPM method of Kofidis and Regalia \cite{KR}) can fail  to find a Z-eigenvalue, see for example, \cite{KR, KM}.  We observed this behavior in  our numerical experiments.    Kolda and Mayo \cite{KM} proved that if  the shift parameter $\alpha<0$  is negative enough (or $\alpha>0$  is large enough),  then $x^{(k)}$ generated by the SSHOPM method converges to a Z-eigenvector. However, the SSHOPM method slows down significantly when a very negative $\alpha<0$ or very large $\alpha>0$  is used \cite{KM}.  Kolda and Mayo \cite{KM} found that using $\alpha=-2, -1, 1,2$ worked well in their tests.  

Since the unconstrained variational principle \reff{Zfun1} leads to negative Z-eigenvalues of $\cA$, we used the SSHOPM method with a negative shift parameter $\alpha$. We solved \reff{Zfun1} via {\tt fminunc} and ran {\tt sshopm} with $\alpha=-2$ on tensors of various dimensions from Example 1 and Example 2. For each tensor we tested, we ran each of {\tt fminunc} and {\tt sshopm} on the tensor 100 times, using a  normalized randomly generated initial vector 
\be
\label{ini_vec}
x_0=\frac{y_0}{\|y_0\|_2},
\ee
where  $y_0={\tt randn}(n,1)$ at each time.   We report the numerical results in Tables 1 and 2, in which  ``CPU time" and ``Accuracy''  denote the ``average CPU time (in seconds)'' and  ``average  $ \hat{\epsilon}=\| \cA \hat{x}^{3} - \lambda \hat{x} \|_2$'' respectively.

\begin{table}
\label{table1}
\footnotesize
\caption{Unconstrained variational principle vs SSHOPM  on some tensors from Example 1 using normalized randomly generated initial vectors} 
\centering 
\begin{tabular}{|c|c|c|c|}\hline
Problem & Method & CPU time &  Accuracy   \\ \hline\hline

 $n=10$ & SSHOPM ($\alpha=-2$) & 0.16  & $5.33 \times 10^{-7} $     \\ \cline{2-4}   
    & \reff{Zfun1} via {\tt fminunc} & 0.12  & $2.12 \times 10^{-9}$  \\ \cline{2-4}    
  \hline 

 $n=20$ & SSHOPM ($\alpha=-2$) & 0.24  &  $4.64 \times 10^{-7}$   \\ \cline{2-4}   
    & \reff{Zfun1} via {\tt fminunc} & 0.18 & $2.65 \times 10^{-9}$     \\ \cline{2-4}    
  \hline 
 
$n=30$ & SSHOPM ($\alpha=-2$) & 0.62  & $3.73 \times 10^{-7}$  \\ \cline{2-4}   
    & \reff{Zfun1} via {\tt fminunc} & 0.34 & $3.70 \times 10^{-9}$ \\ \cline{2-4}    
  \hline 

$n=40$ & SSHOPM ($\alpha=-2$) & 1.69  & $4.05 \times 10^{-7}$    \\ \cline{2-4}   
    & \reff{Zfun1} via {\tt fminunc}& 0.60 & $3.10 \times 10^{-9}$    \\ \cline{2-4}
\hline 

$n=50$ & SSHOPM ($\alpha=-2$) & 4.05  & $3.63 \times 10^{-7}$      \\ \cline{2-4}   
    & \reff{Zfun1} via {\tt fminunc} & 1.36 & $3.30 \times 10^{-9}$   \\ \cline{2-4}    
  \hline 

$n=60$ & SSHOPM ($\alpha=-2$) & 8.24  & $3.08 \times 10^{-7}$     \\ \cline{2-4}   
    & \reff{Zfun1} via {\tt fminunc} & 2.69 & $4.37 \times 10^{-9}$   \\ \cline{2-4}    
  \hline 
\end{tabular}
\end{table}

\begin{table}
\label{table2}
\footnotesize
\caption{Unconstrained variational principle  vs SSHOPM  on some random tensors generated from Example 2 using normalized randomly generated initial vectors} 
\centering 
\begin{tabular}{|c|c|c|c|}\hline
Problem & Method & CPU time &  Accuracy   \\ \hline\hline

 $n=10$ & SSHOPM ($\alpha=-2$) & 0.38  & $1.03 \times 10^{-6} $     \\ \cline{2-4}   
    & \reff{Zfun1} via {\tt fminunc} & 0.18  & $9.98 \times 10^{-9}$  \\ \cline{2-4}    
  \hline 

 $n=20$ & SSHOPM ($\alpha=-2$) & 0.89  &  $1.33 \times 10^{-6}$   \\ \cline{2-4}   
    & \reff{Zfun1} via {\tt fminunc} & 0.27 & $2.19 \times 10^{-8}$     \\ \cline{2-4}    
  \hline 
 
$n=30$ & SSHOPM ($\alpha=-2$) & 2.31  & $1.48 \times 10^{-6}$  \\ \cline{2-4}   
    & \reff{Zfun1} via {\tt fminunc} & 0.53 & $3.65 \times 10^{-8}$ \\ \cline{2-4}    
  \hline 

$n=40$ & SSHOPM ($\alpha=-2$) & 5.87  & $1.61 \times 10^{-6}$    \\ \cline{2-4}   
    & \reff{Zfun1} via {\tt fminunc}& 1.06 & $6.52 \times 10^{-8}$    \\ \cline{2-4}
\hline 

$n=50$ & SSHOPM ($\alpha=-2$) & 11.84  & $1.69 \times 10^{-6}$      \\ \cline{2-4}   
    & \reff{Zfun1} via {\tt fminunc} & 2.21 & $8.88 \times 10^{-8}$   \\ \cline{2-4}    
  \hline 

$n=60$ & SSHOPM ($\alpha=-2$) & 24.28  & $1.79 \times 10^{-6}$     \\ \cline{2-4}   
    & \reff{Zfun1} via {\tt fminunc} & 4.56 & $1.09 \times 10^{-7}$   \\ \cline{2-4}    
  \hline 
\end{tabular}
\end{table}

From Tables 1 and 2, we observe that although both the SSHOPM method and the unconstrained optimization principle \reff{Zfun1} successfully find a Z-eigenvalue of $\cA$ in all cases, solving \reff{Zfun1} via {\tt fminunc} on average uses less CPU time than {\tt sshopm}, particularly when $n$ is large. This is perhaps due to the superlinear convergence property of the BFGS algorithm used in {\tt fminunc}.

A natural question is how to choose a suitable shift parameter $\alpha$ in {\tt sshopm}.  Using $\alpha=-2$ worked well for the problems considered in Tables 1 and 2. However, this choice is not a suitable one for some tensors from Example 3.  We illustrate this in Table 3, in which we report the numerical results on a tensor from Example 3 with dimension $n=25$.  We ran   {\tt sshopm} with various shift parameters and {\tt fminunc}   on this tensor 10 times, using a normalized randomly generated initial vector $x_0$ as defined in \reff{ini_vec} with $n=25$ at each time.   The notations are:
\begin{itemize}
\item Min/Max/Mean Accuracy denotes the minimum, maximum, and mean  $\hat{\epsilon} =  \| \cA \hat{x}^{3} - \lambda \hat{x} \|_2$. 
\item Min/Max/Mean CPU time denotes the minimum, maximum, and mean  CPU time used.
\end{itemize}
Clearly for this example, $\alpha=0,-1,-2,-5,-10, -100$ are not a suitable choice for the shift parameter. Although the tensor used in this test is artificial,  the numerical results  indicate that  choosing a suitable shift parameter can be crucial for the success of the SSHOPM method.

\begin{table}
\label{table3}
\footnotesize
\caption{Unconstrained variational principle vs SSHOPM on a tensor generated from Example 3 with $n=25$, using normalized randomly generated initial vectors} 

\begin{tabular}{|c|c|c|}\hline
 Method & Min/Max/Mean Accuracy &  Min/Max/Mean CPU time \\ \hline \hline
 \reff{Zfun1} via {\tt fminunc} & $1.06 \times 10^{-5} /  3.44 \times 10^{-2} / 4.23 \times 10^{-3} $ & 0.36/1.14/0.63  \\ \hline
SSHOPM ($\alpha=0$) &   $1.36 \times 10^4 /  4.02 \times 10^{4} / 2.91 \times 10^{4} $     &    25.78/26.85/26.09   \\  \hline
SSHOPM ($\alpha=-1$) &   $ 1.08 \times 10^4 / 5.26 \times 10^{4}  / 3.73 \times 10^{4} $     &    25.96/28.58/26.78     \\ \hline
SSHOPM ($\alpha=-2$) &   $ 2.56 \times 10^4 /  5.09 \times 10^{4}/ 3.97 \times 10^{4} $     &    25.84/27.50/26.34      \\ \hline
SSHOPM ($\alpha=-5$) &   $7.30 \times 10^3 /  5.55 \times 10^{4} / 3.58 \times 10^{4} $     &    25.79/27.11/26.25       \\  \hline
SSHOPM ($\alpha=-10$) &    $3.69 \times 10^{3} /  4.55 \times 10^4 / 1.56 \times 10^4 $     &    25.78/25.91/25.81       \\  \hline
SSHOPM ($\alpha=-100$) &    $1.53  /  4.12 \times 10^{4} / 2.80 \times 10^{4} $     &    25.66/ 25.84 /25.75       \\  \hline
SSHOPM ($\alpha=-1000$) &    $6.50 \times 10^{-3} /  1.09 \times 10^{-2} / 8.79 \times 10^{-3} $     &   3.59/4.26/3.81       \\  \hline
\end{tabular}
\end{table}

We now summarize our comparison of the SSHOPM method and the unconstrained variational principle approach:  
\begin{itemize}
\item The SSHOPM method  can be used to find both real and complex Z-eigenvalues and can handle both even and odd order symmetric tensors \cite{KM}.    The unconstrained optimization approach 
can find real $\cB_r$-eigenvalues (including Z-, H-, and D-eigenvalues) and handle  even order (weakly) symmetric tensors. 
\item The unconstrained optimization approach  can be faster than the SSHOPM method particularly when $n$ is large. 
\item If the purpose is to find one Z-eigenvalue for a given tensor, Example 3 shows that choosing a suitable shift parameter may be crucial for the SSHOPM method. The unconstrained optimization approach sometimes needs to solve a shifted problem \reff{shifted1} or \reff{shifted2} (see the next subsection).   However, if tensor $\cA$ has negative (or positive) $Z$-eigenvalues that are well separated from 0,  directly solving Problem \reff{obj1} (or Problem \reff{obj2}) with $\cB=I_{n}^{m/2}$ can generally obtain a negative (or positive) Z-eigenvalue, as Examples 1, 2, and 3 illustrate. 
\item A global polynomial optimization solver can be used to solve the optimization problems arisen in the unconstrained variational principles to find the largest or smallest eigenvalues. The found Z-eigenvalue by the SSHOPM method is not necessarily the largest or the smallest one. 
\end{itemize}

\subsection{Determining positive semidefiniteness}

In some applications, it is important to determine if an even order symmetric tensor $\cA$ is positive semidefinite (see, for example, \cite{QYW, QYX}). 
An attractive property of function \reff{obj1} is that any of its nonzero critical point  is a $\cB_r$-eigenvector of $\cA$ 
corresponding to a $\cB_r$-eigenvalue $\lambda < 0$. This  feature allows us to use a local optimization solver to determine the positive semidefiniteness of an even order symmetric tensor, since $\cA$ is positive semidefinite if all of its H-eigenvalues and all of its Z-eigenvalues are nonnegative ( \cite{Qi}).

We will consider the more general shifted  problem \reff{shifted1} in this subsection. Note that \reff{shifted1} becomes \reff{obj1} when $t=0$. When $t \ne 0$, solving \reff{shifted1} leads to a $\cB_r$-eigenvalue of the tensor $\cA+ t\cB$. In this situation, subtracting $t$ from the computed eigenvalue of  $\cA+t \cB$ will result in a $\cB_r$-eigenvalue of $\cA$.   We consider both $\cB=\cI$ (corresponding to H-eigenvalues) and $\cB=I_{n}^{m/2}$ (corresponding to Z-eigenvalues).

We now summarize the unconstrained optimization approach for determining  the positive semidefiniteness of an even order symmetric tensor $\cA$ in Algorithm 1. 

\newpage
\begin{framed}
\begin{70} 
\label{algorithm1}  
{\rm 
\ \\
\ \\
{\bf Input}: Tensor $\cA$. 
\ \\
{\bf Step 0}. Choose parameters $t < 0$, $0 \leq \eta_1 \leq \eta_2$, and tensor $\cB=\cI$ or $\cB=I_{n}^{m/2}$.      
\ \\
{\bf Step 1}. Solve the unconstrained optimization problem \reff{shifted1} with parameter $t$. Let $\tilde{x}$ and $\tilde{s}$ denote the  optimal solution and optimal objective value of \reff{shifted1} found by the optimization solver.   
\ \\
{\bf Step 2}. 
  \begin{itemize} 
   \item If $ \cB \tilde{x}^{m} \leq \eta_1$, then $\cA$ is positive definite, stop. 
   \item  If $ \cB \tilde{x}^{m} > \eta_2$, then set $\tilde{\lambda} = -  \sqrt{-2m\tilde{s}} - t$ (which is a $\cB_r$-eigenvalue of $\cA$). If $\tilde{\lambda} \geq 0$, then  $\cA$ is positive  semidefinite; otherwise $\cA$ is not positive semidefinite, stop.
\end{itemize}  
}
\end{70}
 \end{framed}

\begin{5} {\rm 
(a) If $\tilde{x}$ is the global minimizer of problem \reff{shifted1} with global minimum value $\tilde{s}$, then Algorithm 1 can always determine the positive semidefiniteness of $\cA$ by choosing $\eta_1=\eta_2=0$. If $ \cB \tilde{x}^{m} \leq \eta_1 =0$, then  
$ \tilde{x} =0$. As $\tilde{x}$ is  the global minimizer  of $s_1$, we have $t+\lambda_{\min} \geq 0$  according to Theorem \ref{mainthm3}. This implies that $\lambda_{\min}   > 0$ since $t < 0$. Thus $\cA$ is positive definite. If $ \cB \tilde{x}^{m} > \eta_2 =0$, then $ \tilde{x} \ne 0$ is an eigenvector corresponding to the eigenvalue $t+\lambda_{\min}<0$. Then  
$\tilde{\lambda} = -\sqrt{-2m\tilde{s}} - t$ gives the minimum eigenvalue $\lambda_{\min}$. Therefore, we can use it to determine the positive 
semidefiniteness of $\cA$.   \\
(b) Ideally, a global polynomial optimization solver should be used  in Step 1. We have found, however, the state of art global polynomial optimization solvers such as {\tt GloptiPoly3} \cite{HLL} cannot handle problem \reff{shifted1} for large $n$ when $m \geq 4$.  On the other hand, a local optimization solver can solve problem \reff{shifted1} for much larger $n$ or $m$. There is no guarantee that Algorithm 1 always successfully determines the positive semidefiniteness of a tensor $\cA$ in this case. As can be seen from our numerical results,  however,  the success rate of using a local optimization solver in Step 1 with a suitable choice of parameter $t$ is quite promising on determining the positive semidefiniteness of even oder symmetric tensors.
}
\end{5}
\begin{5}{\rm 
(a) When Algorithm 1 is implemented, $\eta_1$ and $\eta_2$ are used to numerically determine whether or not 
$\cB \tilde{x}^m$ is $0$. Therefore, $\eta_1$ should be sufficiently small. On the other hand, $\eta_2$ should be a number that is small, but not too small so $ \cB \tilde{x}^{m} > \eta_2$ implies that $\cB \tilde{x}^m \ne 0$. In our numerical experiments, we have found that  $\eta_1=10^{-10}$ and  $\eta_2=10^{-4}$ worked well for our tested examples.      \\ 
(b) If $ \eta_1 <  \cB \tilde{x}^m \leq \eta_2$, then Algorithm 1 is inconclusive. In this case, we can use a different shift parameter $t<0$ and repeat Step 1 and Step 2.  
}
\end{5}

To test the effectiveness of Algorithm 1 for determining the positive semidefiniteness of even order symmetric tensors, we did some numerical experiments in which Problem \reff{shifted1} was solved via {\tt fminunc}.  
The parameters for {\tt fminunc} were the same as in \reff{optset1}. For comparison, we also tested the SSHOPM method.  For this method, we used the same setting as in \reff{optset3} except for that {\tt MaxIts} was changed to ${\tt MaxIts}=10000$.

\begin{60}
\label{example4}
{\rm
The $4^{th}$ order n-dimensional symmetric tensor $\cA$ is generated as follows: First randomly generate tensor 
$ \cT ={\tt randn}(n,n,n,n)$, then  use the {\tt symmetrize} function in the Matlab Tensor Toolbox \cite{BK} to symmetrize $\cT$ and obtain   
$\cZ= {\tt symmetrize} (\cT)$. Finally set 
\bea
 \cA(i,j,k,l) =  \left \{  \begin{array}{cl}
     1000, & {\rm if} \ 1 \leq i=j=k=l \leq n-1; \\
     -1, &   {\rm if} \  i=j=k=l =n; \\
     \cZ(i,j,k,l),  & {\rm otherwise}.
       \end{array}
\right . 
\eea
$\cA$ is not positive semidefinite when $n \geq 2$. 
}
\end{60}

\begin{60} {\rm  The $4^{th}$ order 3-dimensional tensor $\cA$ is defined by
\label{example5}
$\cA(1,1,1,1)=1$; $\cA(2,2,2,2)=0$; $\cA(3,3,3,3)=-0.001$; and $\cA(i,j,k,l)=0$ for all other $i,j,k,l$. This tensor is not positive semidefinite. 
}
\end{60}

We tested Algorithm 1 using  $\cB=\cI$ or $\cB=I_{n}^{m/2}$ and shift parameters $t=0$ or $t=-1$ and compared them with the SSHOPM method with different shift parameters on some tensors from Examples 4 and 5. For each tensor, we ran each of these methods 100 times, using a normalized randomly generated initial vector as defined in \reff{ini_vec} at each time. We used $\eta_1=10^{-10}$ and $\eta_2=10^{-4}$ in Algorithm 1. 

We report the numerical results in Tables 4 and 5. In both tables, ``Success rate'' denotes the percentage of times where a negative eigenvalue was found (and therefore the corresponding method successfully determined that $\cA$ is not positive semidefinite); ``CPU time'' denotes the average CPU time (in seconds); and ``NIT'' denotes the average number of iterations used by the SSHOPM method.   

For the tensor generated from Example 4, the SSHOPM method always converged to the dominate positive eigenvalue when $\alpha=-2$, $\alpha=-10$, $\alpha=-50$. When $\alpha=-100$, it converged to a negative eigenvalue in 34  out of 100 runs. It successfully found a negative eigenvalue when $\alpha=-500$ in all of the 100 runs, using an average CPU time of 30 seconds. On the other hand, Algorithm 1 using $\cB=\cI$ or $\cB=I_{n}^{m/2}$ and $t=0$ or $t=-1$  successfully found a negative eigenvalue in all runs, using much less CPU time.     

For the tensor from Example 5, when $t=0$, Algorithm 1 using $\cB=\cI$ correctly identified that $\cA$ is not positive semidefinite 91\% of times;   Algorithm 1 using $\cB=I_{n}^{m/2}$ was only successful  32\% of times. The failures in both cases were due to that {\tt fminunc} converged to the critical number $x=0$ of \reff{obj1}. This is because {\tt fminunc} is a local optimization solver. It only guarantees to converge to a critical point. We found that using a negative parameter $t$ can significantly increase the success rate, particularly in the $\cB=I_{n}^{m/2}$ case. Indeed, when $t=-1$ was used, Algorithm 1 successfully found a negative eigenvalue 98\% of times in both  $\cB=\cI$ and $\cB=I_{n}^{m/2}$ cases. The two failure runs in each case were due to that the eigenvalue $\lambda=0$ was found.   In all 100 runs, the SSHOPM method (with $\alpha=-2$) terminated when the maximum allowed number of iterations (which is 10000) was reached. In 10 out of 100 runs,  the method terminated at an approximate Z-eigenvalue very close to $\lambda=0$ (and hence failed to correctly determine the positive semidefiniteness of $\cA$).  We plot the computed Z-eigenvalues by the SSHOPM method in the 100 runs in Figure 1. From this figure we observe that the SSHOPM method successfully determined that $\cA$ is not positive semidefinite in less than 90\% of times. 
   
In summary, Algorithm 1 using $\cB=\cI$ or $\cB=I_{n}^{m/2}$ and a negative shift parameter $t$ is more efficient than the SSHOPM method on determining the positive semidefiniteness of the tensors from Examples 4 and 5 we tested in terms of CPU time.  

\begin{table}
\label{table4}
\footnotesize
\caption{Determining positive semidefiniteness using a tensor generated in Example \ref{example4}, $n=30$} 
\centering 
\begin{tabular}{|c|c|c|c|c|}\hline
 Method & Shift parameter & Success rate  &  CPU time &  NIT    \\ \hline \hline

  Algorithm 1 ($\cB=\cI$) & $t=0$ & 100 & 1.36 &     \\ \cline{1-4}  
Algorithm 1 ($\cB=\cI$) & $t=-1$ & 100 & 1.26&     \\ \cline{1-4}  \hline \hline

  Algorithm 1 ($\cB=I_{n}^{m/2}$) & $t=0$  & 100 & 0.99 &     \\ \cline{1-4}
  Algorithm 1  ($\cB=I_{n}^{m/2}$) & $t=-1$ & 100 & 0.99 &     \\ \cline{1-4} \hline \hline

  SSHOPM & $\alpha=-2$              &  0    &  2.97  & 409.4 \\ \cline{1-5} 
  SSHOPM & $\alpha=-10$             &  0    &  2.34  &  310.1  \\ \cline{1-5}   
 SSHOPM &  $\alpha=-50$             &  0    &  7.41  &  1111.2 \\ \cline{1-5} 
 SSHOPM & $\alpha=-100$             &  34   &  3.22  &  449.3 \\ \cline{1-5} 
SSHOPM & $\alpha=-500$              &  100  &  30.00 & 4693.4 \\ \cline{1-5} 
\hline 

\end{tabular}
\end{table}

\begin{table}
\label{table5}
\footnotesize
\caption{Determining positive semidefiniteness using the tensor in Example \ref{example5}} 
\centering 
\begin{tabular}{|c|c|c|c|c|}\hline
 Method & Shift parameter & Success rate &   CPU time &  NIT    \\ \hline \hline

  Algorithm 1 ($\cB=\cI$) & $t=0$ & 91 & 0.28 &     \\ \cline{1-4}  
 Algorithm 1 ($\cB=\cI$) & $t=-1$ & 98 & 0.62 &     \\ \cline{1-4}  \hline \hline

Algorithm 1 ($\cB=I_{n}^{m/2}$) & $t=0$  & 32 & 0.16&     \\ \cline{1-4}
Algorithm 1  ($\cB=I_{n}^{m/2}$) & $t=-1$ & 98 & 0.47 &     \\ \cline{1-4} \hline \hline

  SSHOPM & $\alpha=-2$           &      $\leq 90$  &  34.58  & 10000 \\ \cline{1-5} 
\hline 

\end{tabular}
\end{table}

\begin{figure}[!h]
\label{figure1}
\begin{center}
\includegraphics[scale=.35]{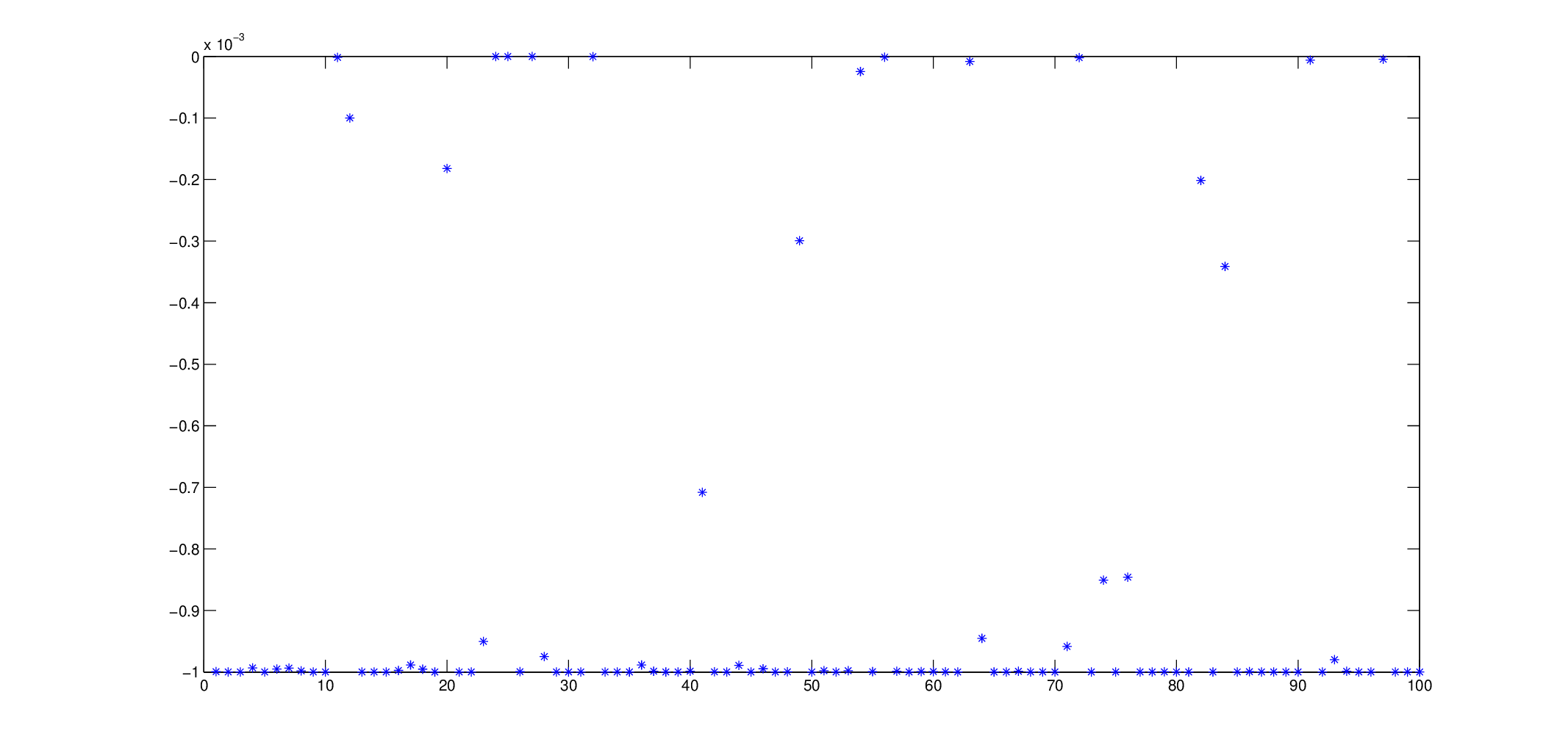}
\caption{The computed Z-eigenvalues by the SSHOPM method in the 100 runs on tensor $\cA$ from Example 5.} 
\end{center}
\end{figure}

We now turn to the performance of Algorithm 1 on positive semidefinite (definite) tensors. We tested some tensors from Examples \ref{example6} and \ref{example7} using Algorithm 1 with  $\cB=\cI$ and $t=-1$ and  with  $\cB=I_{n}^{m/2}$ and $t=-1$. For each tensor we tested, we run each method 100 times, using a normalized randomly generated initial vector as defined in \reff{ini_vec} at each time. We report the numerical results in Tables 6 and 7. In Table 6, ``min $\cB \tilde{x}^m$'' denotes the smallest $\cB \tilde{x}^m$ in 100 runs and ``Success rate'' denotes the percentage of times where the minimum eigenvalue $\lambda=0$ was obtained. In Table 7,  ``min/max $ \cB \tilde{x}^m $'' denotes the smallest and largest $ \cB \tilde{x}^m $ in 100 runs and  ``Success rate'' denotes the percentage of times when the method correctly determined that $\cA$ is positive definite. From Tables 6 and 7, we observe that Algorithm 1 using  ``$\cB=\cI$ and $t=-1$'' or   ``$\cB=I_{n}^{m/2}$ and $t=-1$'' was able to efficiently determine the positive semidefiniteness of the tensors from Examples 6 and 7 we tested.

\begin{60} {\rm  The $4^{th}$ order n-dimensional tensor $\cA$ is defined by
\label{example6}
$\cA(k,k,k,k)={\tt rand}(1)$ for $k=1,2,\cdots, n-1$; $\cA(n,n,n,n)=0$;  and $\cA(i,j,k,l)=0$ for all other $i,j,k,l$. 
$\cA$ is positive semidefinite. }
\end{60}

\begin{60} {\rm  Consider the positive definite $4^{th}$ order n-dimensional tensor $\cA$  defined by
\label{example7}
$\cA(k,k,k,k)=10k$ for $k=1,2,\cdots, n$;   and $\cA(i,j,k,l)=0$ for all other $i,j,k,l$. 
}
\end{60}

\begin{table}[!h]
\label{table6}
\footnotesize
\caption{Determining positive semidefiniteness using a tensor in Example \ref{example6}, $n=30$} 
\centering 
\begin{tabular}{|c|c|c|c|c|}\hline
 Method                &   min $\cB \tilde{x}^m $ & Success rate    &  CPU time     \\ \hline 
  Algorithm 1 ($\cB=\cI$, $t=-1$) &           1.00           & 100  &   1.68  \\ \cline{1-4}  \hline 
  Algorithm 1  ($\cB=I_{n}^{m/2}$, $t=-1$) & 1.00    &   100   &  2.54   \\ \cline{1-4} \hline 
\end{tabular}
\end{table}

\begin{table}[!h]
\label{table6}
\footnotesize
\caption{Determining positive semidefiniteness using a tensor in Example \ref{example7}, $n=30$} 
\centering 
\begin{tabular}{|c|c|c|c|c|}\hline
 Method  &  min/max $ \cB \tilde{x}^m $ &  Success rate   &  CPU time     \\ \hline 
  Algorithm 1 ($\cB=\cI$, $t=-1$) &   $7.12 \times 10^{-17} $ / $ 2.65 \times 10^{-15} $     &    100   &   1.66   \\ \cline{1-4}  \hline 
  Algorithm 1  ($\cB=I_{n}^{m/2}$, $t=-1$) & $4.98 \times 10^{-17} $ / $ 1.97 \times 10^{-13} $ &  100    &  1.81   \\ \cline{1-4} \hline 
\end{tabular}
\end{table}

\section{Final Remarks}

We have introduced  two unconstrained optimization problems and obtained some variational characterizations for the minimum and maximum $\cB_r$ 
 eigenvalues of an even order weakly symmetric tensor, where $\cB$ is weakly symmetric positive definite. These unconstrained optimization problems can be solved using some powerful optimization algorithms, such as the BFGS method. This approach can be used to find a Z-,  H-, and D-eigenvalue  of an even order weakly symmetric tensor.  We have provided some numerical results indicating that our approach of solving Problem \reff{Zfun1} via {\tt fminunc} compares favorably to the approach of solving \reff{minconopt} via {\tt fmincon} and the SSHOPM method for finding a Z-eigenvalue of an even order symmetric tensor.  Furthermore, we have  provided some  numerical results that show  the unconstrained optimization approach is promising on determining positive semidefiniteness of an even order symmetric tensor.

A direction for future research is to develop a global polynomial optimization algorithm that can solve problems \reff{obj1} and \reff{obj2} and their shifted versions when  $n$ (and/or $m$) is large. \\

\ \\
{\bf Acknowledgment.} The author is very grateful to the referees for their constructive comments and suggestions, which have helped improve the content and presentation of the paper. 



\begin{thebibliography}{999}


\bibitem{Au89} G. Auchmuty, Unconstrained variational principles for eigenvalues of real symmetric matrices, {\it SIAM Journal on Mathematical Analysis}, 1989, 20(5): 1186--1207.

\bibitem{Au91} G. Auchmuty, Globally and rapidly convergent algorithms for  symmetric eigenproblems, {\it SIAM Journal on Matrix Analysis and Applications}, 1991, 12(4): 690--706.

\bibitem{BK} B.W. Bader, T.G. Kolda and others, MATLAB Tensor Toolbox Version 2.5, 2012. URL: \\ http://www.sandia.gov/~tgkolda/TensorToolbox/  

\bibitem{CS} D. Cartwright and B. Sturmfels, The number of eigenvalues of a tensor, {\it Linear Algebra and its Applications}, 2013, 438(2): 942--952. 

\bibitem{CPZ08} K.C. Chang, K. Pearson and T. Zhang, Perron-Frobenius theorem for nonnegative tensors, {\it Commun. Math. Sci.}, 2008, 6(5): 507--520.


\bibitem{CPZ09} K.C. Chang, K. Pearson and T. Zhang, On eigenvalues of real symmetric tensors, {\it J. Math. Anal. Appl.}, 2009, 350: 416--422.

\bibitem{Dai} Y. Dai and C. Hao, A subspace projection method for finding the extreme Z-eigenvalues of supersymmetric positive definite tensor, a talk given at the International Conference on the Spectral Theory of Tensors, Nankai University, 2012. 

\bibitem{FGH} S. Friedland, S. Gaubert and L. Han, Perron-Frobenius theorem for nonnegative multilinear forms and extensions, {\it Linear Algebra and Applications}, 2013, 438(2): 738--749.

\bibitem{HLL} D. Henrion, J.-B. Lasserre, and J. L\"{o}fberg, GloptiPoly3: moments, optimization and semidefinite programming, {\it Optimization Methods and Software}, 2009, 24: 761--779.


\bibitem{KR} E. Kofidis and Ph. Regalia,  On the best rank-1 approximation of higher-order supersymmetric tensors, {\it SIAM J. Matrix Anal. Appl.}, 2002, 23:  863--884. 

\bibitem{KM} T.G. Kolda and J.R. Mayo, Shifted power method for computing tensor eigenpairs, {\it SIAM J. Matrix Analysis and Applications}, 2011, 32: 1095--1124.  

\bibitem{LQY} G. Li, L. Qi, and G. Yu, The Z-eigenvalues of a symmetric tensor and its application to spectral hypergraph theory, Department of Applied Mathematics, University of
New South Wales, December 2011.

\bibitem{Lim} L.-H. Lim, Singular values and eigenvalues of tensors: a variational approach, {\it Proceedings of the IEEE International Workshop on Computational Advances
in Multi-Sensor Adaptive Processing (CAMSAP'05)}, 2005, 1: 129--132.

\bibitem{Mathworks} The Mathworks, Matlab  7.8.0, 2009. 

\bibitem{NW} J. Nocedal and S. Wright, {\it Numerical Optimization}, 2nd ed., Springer, 2006.

\bibitem{PSU} A. L. Peressini, F. E. Sullivan, J. J. Uhl, {\it The Mathematics of Nonlinear Programming}, Springer, 1988.  


\bibitem{Qi} L. Qi, Eigenvalues of a real supersymmetric tensor, {\it J. Symb. Comput.}, 2005, 40:
1302--1324.

\bibitem{QSW} L. Qi, W. Sun, and Y. Wang, Numerical multilinear algebra and its applications, {\it Frontiers of Mathematics in China}, 2007, 2: 501--526.

\bibitem{QWWang} L. Qi, F. Wang, and Y. Wang, Z-eigenvalue methods for a global optimization polynomial optimization problem, {\it Mathematical Programming}, 2009, 118: 301--306. 

\bibitem{QWWu} L. Qi, Y. Wang, and E.X. Wu, D-eigenvalues of diffusion kurtosis tensors, {\it Journal of Computational and Applied Mathematics}, 2008, 221: 150--157.  


\bibitem{QYW} L. Qi, G. Yu, and E.X. Wu, Higher order positive semi-definite diffusion tensor imaging, {\it SIAM Journal on Imaging Sciences}, 2010, 3: 416--433. 

\bibitem{QYX} L. Qi, G. Yu, and Y. Xu, Nonnegative diffusion orientation distribution function, {\it Journal of Mathematical Imaging and Vision}, 2013, 45(2): 103--113. 

\end{thebibliography}
\end{document}